# Rate-induced tipping in a solvable model with the Allee effect

Hidekazu Yoshioka[1, *]

[1]Japan Advanced Institute of Science and Technology, 1-1 Asahidai, Nomi, Ishikawa 923-1292, Japan. ORCID: 0000-0002-5293-3246

* Corresponding author: yoshih@jaist.ac.jp

***Abstract***

We present a novel exactly solvable ordinary differential equation model for rate-induced tipping: a dynamic phenomenon of dynamical systems where a time-dependent parameter triggers the transition of stability of a system. Our model contains an Allee effect that induces a saddle point and admits an explicit solution along with the extinction threshold of a time-dependent Allee parameter. More specifically, we derive an integral inequality that serves as a necessary condition for the occurrence of rate-induced tipping. A remarkable point in the proposed model is that it can handle population extinction such that the solution completely vanishes in a finite amount of time. An unconditionally stable cubature method suitable for our model is proposed, and its superiority over the classical forward Euler method is discussed. As a proof of concept, the model reproduces the rise and fall of the membership of inland fisheries cooperatives in Japan from 1963 to 2023. The proposed model serves as a tractable tool for studying R-tipping phenomena.

***Keywords***

Allee effect; Rate-induced tipping; Extinction time; Application to the fisheries sector

## 1. Introduction

### 1.1 Background

Rate-induced tipping (R-tipping) is a dynamic phenomenon triggered by a parameter change in a dynamic system where the stability of the system varies over time [1]. R-tipping has been studied in detail for natural and human systems because humans face this phenomenon in various parts of the world [2,3]. The stability of the dynamics of the Earth system is one such example where the research interests are large-scale phenomena such as the Atlantic Meridional Overturning Circulation (AMOC) [4] and tropical cyclones [5,6]. R-tipping in the Lorenz model, a representative model for climate dynamics, has been evaluated extensively [7]. Cascade tipping phenomena between the Greenland ice sheet and AMOC have been extensively studied with possible tipping scenarios [8]. A conceptual model for tracking the coevolution between life and the environment and their collapse has been proposed [9]. Local environmental phenomena, such as zombie fires [10] and the soil temperature increases (called bomb instability) [11], have also been studied. Rate-induced transitions in view of agroforestry under rapidly changing climates have been extensively reviewed [12].

R-tipping has also been investigated in the context of ecology, where the study interests are the coexistence and extinction of species. The tipping stability of two-patch population dynamics was investigated with a focus on regime shifts [13]. R-tipping in a deterministic population dynamics model of complex networks has been studied by focusing on the species extinction rate [14]. R-tipping focusing on biological control for species conservation has been analyzed with an emphasis on the cost and benefit of control and its rate [15]. R-tipping phenomena has also been studied in other contexts, such as synchronization [16], multiagent systems with cooperation [17], musical applications [18], and deep learning methods for tipping prediction [19,20].

A nominal model for studying R-tipping of continuous-time dynamical systems is the following Allee model (nominal model) for describing population dynamics with a saddle point:

$$\underbrace{\frac{\mathrm{d}X_t}{\mathrm{d}t}}_{\text{Speed of change}} = \underbrace{rX_t}_{\text{Proportional growth}} \underbrace{\left(K - X_t\right)}_{\text{Growth saturation}} \underbrace{\left(X_t - a\right)}_{\text{Allee efect}}, \quad t > 0 \tag{1}$$

subject to an initial condition $X_0 \in [0, K]$ with parameters $r > 0$ and $0 < a < K$ (e.g., Eq. (1) in Feudel [21]; Eq. (7.1) Krakovská et al. [22]). The solution $X = \left(X_t\right)_{t\geq 0}$ to the ODE (1) typically represents some population in a habitat. The stationary equilibria $X_\infty = \lim_{t\to+\infty} X_t$ of (1) are $X_\infty = 0, K$ (stable) and $X_\infty = a$ (unstable), and any solutions converge to $X_\infty = 0$ if $X_0 \in [0, a)$ and to $X_\infty = K$ if $X_0 \in (a, K]$. In the context of R-tipping, the Allee parameter $a$ is replaced by a time-dependent coefficient. In this case, the limit $X_\infty$ depends not only on the initial condition but also on the Allee parameter, and the dynamics of the ODE (1) become more complicated; in particular, time changes of $a$ control the destination of the solution $X$ in the long run, where the stability change occurs when $X$ and $a$ intersect [21].

Although the ODE (1) is fairly simple with explicit stationary equilibria, its drawback is that it does not admit explicit solutions even for constant $a$ cases, which is partly due to the cubic nonlinearity on the right-hand side. Indeed, with a quadratic right-hand side, similar ODEs (i.e., logistic equations) can be solved exactly. Another characteristic, which is not necessarily a drawback, of (1) is that its equilibria are not attained in a finite time but are only the limit for $t \to +\infty$. This implies that the ODE (1) is not applicable to problems where the extinction time is of interest [23]. These limitations in the nominal ODE (1) motivated us to explore a new ODE model that has a stability property close to (1) but with an explicit solution and the ability to address finite-time extinctions.

### 1.2 Aim and contribution

The aim of this paper is to propose and analytically and computationally investigate an ODE model with an Allee effect. **Table 1** qualitatively compares the two ODEs (1) (nominal model) and (2) (proposed model explained in detail in **Section 2**); the key difference between them is the finiteness of the extinction time (smallest $t$ such that $X_t = 0$ when $X_0 > 0$) and the availability of explicit solutions. We show that the explicit solution to the proposed ODE is found through a variable transformation method, and the extinction time is found as a byproduct of the derivation procedure of the solution. Moreover, we derive a threshold condition to judge whether some R-tipping occurs depending on the initial conditions and model parameters. To the best of the author's knowledge, existing ODEs with Allee effects in the context of R-tipping generally do not have explicit solutions. Our ODE therefore serves as a novel analytical approach for studying R-tipping phenomena.

Owing to including the logarithm of the solution in the coefficient, the proposed model is also related to, and actually was motivated from Gompertz-type ODEs, which are major models for describing population dynamics that often admit closed-form solutions [24-27]. A variety of Gompertz-type ODEs have been studied thus far, but their possible connection to R-tipping phenomena has yet to be well studied. We can exploit the appearance of logarithms of the solution in the coefficient, as for Gompertz-type ODEs, and have access to the explicit formula of the proposed model. More specifically, a key technical point in the proposed ODE is that it can be reformulated as a logistic-type ODE by a suitable transformation of variables, which can be further reduced to a linear ODE by another transformation. These transformations apply to both constant and time-dependent Allee parameter cases. We can therefore study R-tipping phenomena in a tractable way.

We propose a cubature method that directly computes some integrals in the exact solution. This numerical method is time explicit and unconditionally stable (always preserves the nonnegativity of solutions) and is more advantageous than classical first-order numerical methods for ODEs, such as the forward Euler method, because it does not depend on the regularity conditions of ODE coefficients. We show that the logarithmic singularity of the coefficients in the proposed ODE at the boundary of the domain degrades the computational accuracy of the extinction time in the forward Euler method.

We also address an application of the proposed model to a problem in inland fisheries, where we investigate time series data of membership in inland fisheries in Japan, a major industry sector contributing to the economy, environment, and ecosystem of the country. The time series data have increasing and decreasing trends, and we study them through the lens of R-tipping, providing a new interpretation of the data. Consequently, this paper contributes to the formulation, analysis, and application of a new mathematical model for R-tipping phenomena.

**Table 1.** Qualitative difference and similarity between the two ODEs (1) and (2).

| | ODE (1) (nominal model) | ODE (2) (proposed model) |
|---|---|---|
| Stable equilibria (constant $a$ case) | $X_\infty = 0, K$ | |
| Unstable equilibria (constant $a$ case) | $X_\infty = a$ | |
| Right-hand side | Cubic $rX_t(K - X_t)(X_t - a)$ | Nonpolynomial $rX_t(\ln K - \ln X_t)(\ln X_t - \ln a)$ |
| Explicit solution | Not available | Available even for time-dependent $a$ |
| Finite-time extinction | No | Yes |

### 1.3 Structure of this paper

**Section 2** presents our ODE model and derives its solution. This section also studies the R-tipping and finite-time extinction of the solution. **Section 3** addresses computational investigations, including an application of the proposed model. **Section 4** summarizes our obtained results and presents future perspectives on this study. **Appendix** presents the auxiliary results.

## 2. ODE model

### 2.1 Formulation

Let $t \geq 0$ be time. In the sequel, we use the following convention: $0(\ln 0)^n = 0$ ($n = 1, 2$). The ODE model proposed in this paper is as follows:

$$\underbrace{\frac{\mathrm{d}X_t}{\mathrm{d}t}}_{\text{Speed of change}} = \underbrace{rX_t}_{\text{Proportional growth}} \underbrace{(\ln K - \ln X_t)}_{\text{Growth saturation}} \underbrace{(\ln X_t - \ln a)}_{\text{Allee efect}}, \quad t > 0 \tag{2}$$

subject to an initial condition $X_0 \in [0, K]$ with parameters $r > 0$ and $0 < a < K$ as for the ODE (1). We assume that $r$ and $K$ are constants but that $a$ is possibly a time-dependent coefficient. If $X_t$ represents some population in a habitat, then the parameter $r$ is its growth rate, $K$ is the carrying capacity, and $a$ is the Allee parameter (saddle parameter) that separates behaviors $X_t$ in the long run. If $a$ is a constant, then the stationary equilibria $X_\infty$ of (2) are $X_\infty = 0, K$ (stable) and $X_\infty = a$ (unstable), and for (1), any solutions to (2) converge to $X_\infty = 0$ if $X_0 \in [0, a)$ and to $X_\infty = K$ if $X_0 \in (a, K]$. Therefore, the two ODEs (1) and (2) share qualitatively the same stability property (see also **Table 1**). The proposed model is reduced to a Gompertz-type ODE if the factor $\ln X_t - \ln a$ is omitted.

In the rest of this paper, we assume that the Allee parameter $a$ is time-dependent unless otherwise specified because we are interested in the nonstationary dynamics of the ODE (2), and we explicitly indicate the time dependence of the Allee parameter as $a_t$. We assume that $a_t$ is continuous and $0 < a_t < K$ for all $t > 0$.

### 2.2 Explicit solution

We derive the explicit solution to the proposed ODE (2). For later use, we set the following quantity: if $X_0 \in (0, K]$,

$$I_t = \frac{1}{\ln(K) - \ln(X_0)} - r \underbrace{\int_0^t e^{-r\int_0^s \ln\left(\frac{K}{a_u}\right) \mathrm{d}u} \mathrm{d}s}_{L_t}, \quad t > 0. \tag{3}$$

We have $I_0 > 0$ and understand that $I_t = +\infty$ if $X_0 = K$. We also set

$$\tau = \inf\{t > 0 : I_t = 0\}. \tag{4}$$

We understand that $\tau = +\infty$ if $I_t > 0$ for all $t > 0$. The integral in (3), which is called $L_t$ in the sequel, is increasing for all $t > 0$, and hence $I_t$ is decreasing for all $t > 0$.

The following proposition is the main theoretical result in this paper.

***Proposition 1***

*If $X_0 \in (0, K]$, then the ODE (2) admits the following unique solution:*

$$X_t = \begin{cases} K \exp\left(-\dfrac{1}{I_t} e^{-r\int_0^t \ln\left(\frac{K}{a_s}\right) \mathrm{d}s}\right) & (0 \le t < \tau) \\ 0 & (t \ge \tau) \end{cases}, \quad t \ge 0 \tag{5}$$

*and the extinction time of this solution is given by $\tau$ in (4). If $X_0 = 0$, then the unique solution is $X_t = 0$ $(t \ge 0)$. Moreover, the following ordering property with respect to the initial condition holds true:*

$$X_t\big|_{X_0 = x_0} \le X_t\big|_{X_0 = x_1}, \quad 0 \le x_0 \le x_1 \le K, \quad t > 0. \tag{6}$$

Its proof is placed in **Appendix**. Its strategy is lengthy but simple, where we transform the ODE (2) to a linear one by using the new variable $W_t = (\ln K - \ln X_t)^{-1}$. Inspecting the proof shows that the extinction time $\tau$ is exactly the time at which $I_t$ changes the sign from positive to negative.

The second part of **Proposition 1** shows that the occurrence of extinction of $X$ and its timing can be determined by $I_t$ because it does not depend on $X$ itself. Another strength of the result obtained is that we do not have to assume the functional form of $a_t$, suggesting that it applies to a wide variety of time dependences of $a_t$ as we study in **Section 3**.

An important result drawn from **Proposition 1** is that a finite-time extinction occurs if and only if the cumulative impact of the Allee effect, which is evaluated by $L_\infty$, is larger than a function of the initial condition $X_0$. This means that what is critical for determining the extinction is not instantaneous values of $a_t$ but rather their accumulation. Moreover, **Proposition 1** shows that $\tau < +\infty$ ($I_t \le 0$ at some $t > 0$) is a necessary condition for the R-tipping where $X_t$ crosses $a_t$ from top to down at least once and $X_\infty = 0$. If in addition $a_t$ is strictly increasing for all $t > 0$ and $X_0 > a_0$, then $I_t \le 0$ at some $t > 0$ is the necessary and sufficient condition for the occurrence of R-tipping.

Finally, the finite-time extinction of the solution to the ODE (2) is caused by the non-Lipschitz nature of its right-hand side; indeed, an elementary calculation shows that

$$\frac{\mathrm{d}}{\mathrm{d}X_t}\left(X_t\left(\ln K - \ln X_t\right)\left(\ln X_t - \ln a\right)\right) = O\left(\left(\ln X_t\right)^2\right) \text{ if } X_t > 0 \text{ is small,} \tag{7}$$

which highlights a qualitative difference between the nominal and proposed models because the former has a right-hand side that is Lipschitz continuous for all $0 \le X_t \le K$.

***Remark 1.*** For the constant $a$ case, we obtain the extinction time

$$\tau = \frac{1}{r\left(\ln(K) - \ln(a)\right)} \ln\left(\frac{\ln(K) - \ln(X_0)}{\ln(a) - \ln(X_0)}\right) \text{ if } X_0 \in (0, a). \tag{8}$$

The solution for the constant $a$ case is fully explicit because we can analytically obtain all the integrals in (5) and (3): for $0 \le t < \tau$,

$$X_t = K\exp\left(-e^{-r(\ln(K)-\ln(a))t}\left\{\left(\ln(K) - \ln(X_0)\right)^{-1} - \left(\ln(K) - \ln(a)\right)^{-1}\left(1 - e^{-r(\ln(K)-\ln(a))t}\right)\right\}^{-1}\right). \tag{9}$$

We use this exact solution to validate our numerical methods.

## 3. Investigations

In this section, first, we study the solution behavior and its extinction time of the proposed model by assuming sigmoidal and oscillatory profiles. Second, we address an industrial application of the model.

### 3.1 Numerical methods

Integrals appearing in (5) are not always explicitly available. We therefore propose the use of two numerical methods. The first method is the classical time-explicit Euler method, which directly discretizes the ODE (2) (called the Euler method in the sequel):

$$X_{k+1} = X_k + hrX_k\left(\ln K - \ln X_k\right)\left(\ln X_k - \ln a_{t_{k+1}}\right), \quad k = 0, 1, 2, \ldots, \tag{10}$$

if $X_k > 0$, where $h > 0$ is the time step and $X_k$ represents the numerical solution at time $t_k = kh$. We set $X_k = 0$ if $k \ge k'$ where $k' = \inf\left(k \in \mathbb{N} : X_k < 0\right)$. Then, we obtain the numerical $\tau$ as $k'h$. As it is well-known, we must choose a sufficiently small $h$ for computational stability (e.g., Chapter 11.3.3 in Quarteroni et al. [28]). This is problematic for the proposed model (2) because of its non-Lipschitz nature of the right-hand side at $X_t = 0$. At the same time, the accuracy is also expected to degrade near $X_t = 0$. This point will be discussed in our computational investigations of extinction times.

The second numerical method directly evaluates $L_t$, the integral term in (3), via a direct quadrature. This method (called the quadrature method in the sequel) is based on

$$\left.\frac{\mathrm{d}L_t}{\mathrm{d}t}\right|_{t=t_k} = e^{-r\int_0^{t_k}\ln\left(\frac{K}{a_u}\right)\mathrm{d}u} \approx \exp\left(-rh\sum_{j=1}^{k}\frac{1}{2}\left(\ln\left(\frac{K}{a_{t_{j-1}}}\right)+\ln\left(\frac{K}{a_{t_j}}\right)\right)\right),\quad k = 1,2,3,\ldots. \tag{11}$$

We obtain the approximation $I_k$ of $I_{t_k}$ as follows:

$$I_k = \frac{1}{\ln K - \ln X_0} - rh\sum_{j=1}^{k}\left.\frac{\mathrm{d}L_t}{\mathrm{d}t}\right|_{t=t_j},\quad k = 1,2,3,\ldots. \tag{12}$$

Then, we obtain numerical $\tau$ as $k''h$ where $k'' = \inf\left(k \in \mathbb{N} : I_k < 0\right)$, and finally $X_k$ as

$$X_k = K\exp\left(-\left.\frac{\mathrm{d}L_t}{\mathrm{d}t}\right|_{t=t_k}\frac{1}{I_k}\right),\quad k = 1,2,3,\ldots \tag{13}$$

if $I_k > 0$ and $X_k = 0$ if $I_k < 0$. The cubature method is essentially a time-explicit discretization of the formal ODE of the variable $L_t$, as indicated in (12).

The following proposition shows that numerical solutions computed by using the cubature method is always nonnegative irrespective of $h > 0$.

***Proposition 2***

*With the cubature method,* $X_k \geq 0$ *for all* $k = 1,2,3,\ldots$.

Its proof is placed in **Appendix**. This proposition highlights an advantage of the cubature method. In contrast, a drawback of the cubature method is that it specializes in the proposed model and is not applicable to the other ODE models. In this view, the Euler method is more versatile.

We report computational results about a constant $a\left(=0.5\right)$ case with $X_0 = 0.32$, $r = 1$, and $K = 1$. **Table A1 in Appendix** compares the accuracy of the Euler and cubature methods against the extinction time $\tau$ for a constant $a$ case and shows that the cubature method (approximately first-order accuracy) performs better than the Euler method (approximately 0.5th-order accuracy), demonstrating the superior accuracy of the cubature method. **Table A2** for errors in the solution shows that each method has first-order accuracy and the error is slightly larger in the cubature method at the chosen initial condition; however, we also find that the size of the errors is reversed for a smaller $X_0$ such as $X_0 = 0.16$ but remains comparable (**Table A4**); however, this reversal is not observed for $\tau$ (**Table A3**). These observations suggest that both the Euler and cubature methods have comparable accuracies for the solution, but the extinction time $\tau$ can be estimated more accurately via the cubature method. These results are considered due to the regularity of the right-hand side of (2); it is only locally Lipschitz continuous for $X_t > 0$ but has a slope with the logarithmic singularity $O\left(\left(\ln X_t\right)^2\right)$ near $X_t = 0$.

### 3.2 Sigmoidal case

We investigate the ODE (2) and the hitting time $\tau$ in (4), where the Allee parameter $a_t$ is given by the sigmoid function (e.g., Ritchie and Sieber [29]; Slyman and Jones [30]):

$$a_t = \frac{\bar{a} - \underline{a}}{1 + \exp\left(\pm \frac{t-\theta}{\varepsilon}\right)} + \underline{a}, \quad t > 0 \tag{14}$$

where $\bar{a}$ and $\underline{a}$ are the upper and lower bounds of $a_t$, $\theta \in \mathbb{R}$ is a shift parameter, and $\varepsilon$ is the width parameter. The positive and negative signs in the exponential function in (14) correspond to decreasing and increasing $a_t$ in time, respectively. The parameters $\bar{a}$ and $\underline{a}$ control the size of $a_t$, $\theta$ is the time of transition, and $\varepsilon$ is the width of the transition. For illustration purposes, we set the following parameter values: $r = 1$, $K = 1$, $\theta = 1$, $\bar{a} = 0.9$, $\underline{a} = 0.1$, and $\varepsilon = 0.1$, with which we can generate a variety of trajectories converging to either $X_\infty = K$ or $X_\infty = 0$ depending on the initial condition $X_0$.

**Figure 1** shows the computed trajectories of $X_t$ for initial conditions $X_0 = 0.04i$ ($i = 0,1,2,...,25$) via the cubature method with $h = 0.0001$. We do not present the results of the Euler method because they are visually the same. **Figure 1(a)** for the increasing $a_t$ case shows that the trajectories of $X_t$ above that of $a_t$ converge to the stationary equilibrium $X_\infty = K$, whereas the others converge to $X_\infty = 0$. This R-tipping phenomenon is due to the increasing nature of $a_t$, where the crossing of $X_t$ with $a_t$ means that we have $X_t < a_t$ and hence $\frac{\mathrm{d}X_t}{\mathrm{d}t} < 0$ strictly after the crossing time, implying convergence to $X_\infty = 0$. A converse phenomenon is observed in **Figure 1(b)** for the decreasing $a_t$ case where the crossing of $X_t$ with $a_t$ means that we have $X_t > a_t$ and hence $\frac{\mathrm{d}X_t}{\mathrm{d}t} > 0$ strictly after the crossing time, implying convergence to $X_\infty = K$. In both cases, trajectories $X_t$ that cross (resp., that do not cross) $a_t$ are nonmonotone (resp., monotone) in time.

Next, we compare the computed $\tau$ between the Euler and quadrature methods. **Table 2** shows the computed $\tau$ corresponding to the cases depicted in **Figure 1**, demonstrating that there is a reasonable agreement between the two methods. Employing a finer resolution $h = 0.00001$ yields analogous results as demonstrated in **Table 3**. In both **Tables 2 and 3**, the Euler method yields a smaller $\tau$ than the cubature method does. By **Tables A1 and A3 in Appendix**, we can assume that the cubature method more accurately computes $\tau$ than does the Euler method for both the increasing and decreasing $a_t$ cases as well; indeed, the computed $\tau$ changes by the cubature method at most 0.03% by improving the computational resolution, whereas those computed by the cubature method change by 0.1% to 3%. This observation suggests faster convergence of the former.

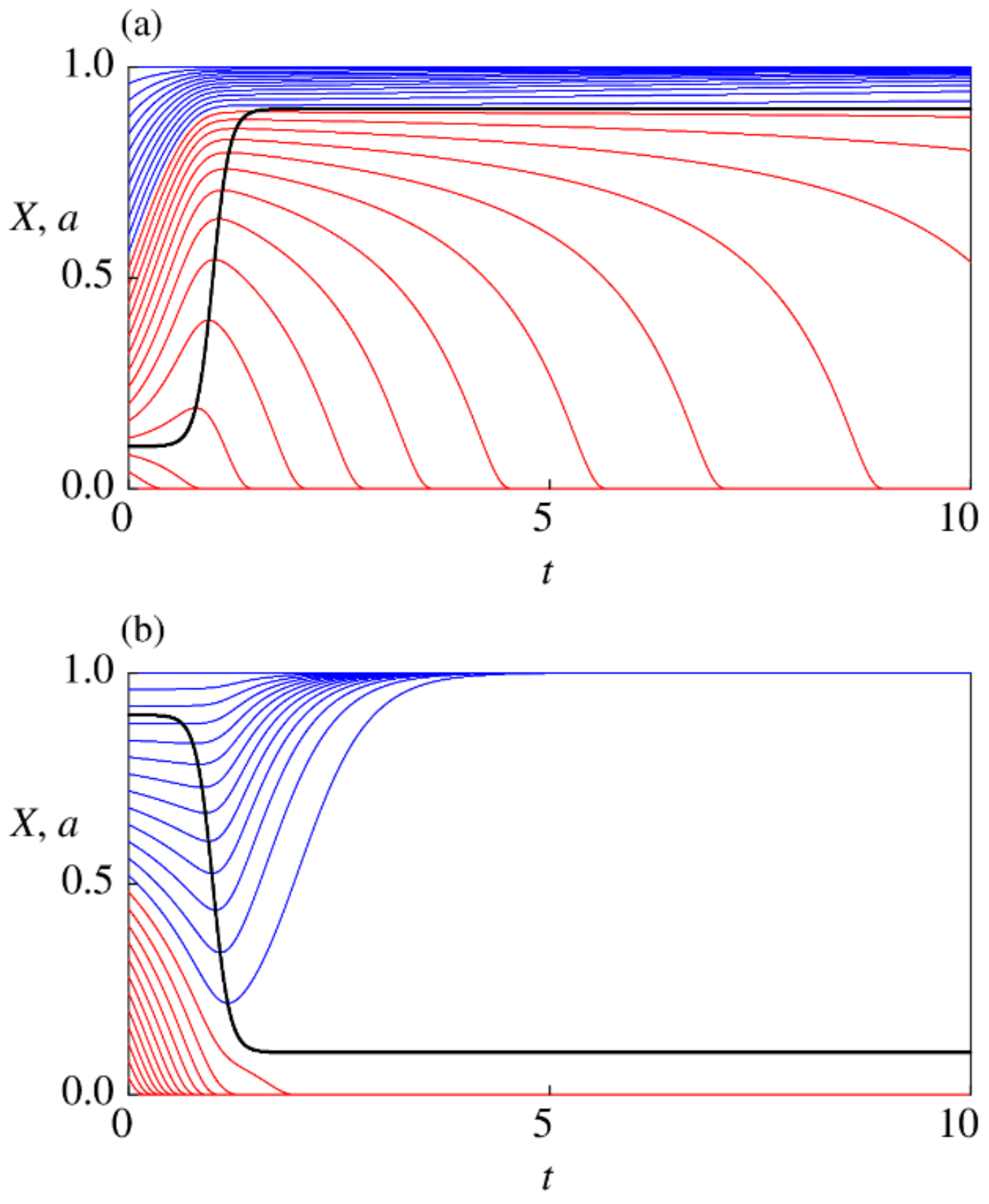

**Figure 1.** Computed $X_t$ for the (a) increasing and (b) decreasing $a_t$ cases via the cubature method. The blue and red curves show the trajectories of $X_t$ with $X_\infty = 1$ (blue) and $X_\infty = 0$ (red). The black curve in each figure panel represents the corresponding trajectory of $a_t$ .

**Table 2.** Comparison of the computed $\tau$ between the Euler and cubature methods: $h = 0.0001$.

| $X_0$ | Increasing $a_t$ case | | Decreasing $a_t$ case | |
|---|---|---|---|---|
| | Euler | Cubature | Euler | Cubature |
| 0.04 | 0.5299 | 0.5448 | 0.3011 | 0.3159 |
| 0.08 | 0.9721 | 0.9871 | 0.3897 | 0.4045 |
| 0.12 | 1.5756 | 1.5910 | 0.4690 | 0.4839 |
| 0.16 | 2.2155 | 2.2310 | 0.5473 | 0.5622 |
| 0.20 | 2.9174 | 2.9330 | 0.6281 | 0.6429 |
| 0.24 | 3.7142 | 3.7300 | 0.7137 | 0.7286 |
| 0.28 | 4.6472 | 4.6631 | 0.8067 | 0.8215 |
| 0.32 | 5.7753 | 5.7914 | 0.9101 | 0.9249 |
| 0.36 | 7.1927 | 7.2090 | 1.0298 | 1.0447 |
| 0.40 | 9.0666 | 9.0831 | 1.1812 | 1.1960 |

**Table 3.** Comparison of the computed $\tau$ between the Euler and cubature methods: $h = 0.00001$.

| $X_0$ | Increasing $a_t$ case | | Decreasing $a_t$ case | |
|---|---|---|---|---|
| | Euler | Cubature | Euler | Cubature |
| 0.04 | 0.53983 | 0.54466 | 0.31106 | 0.31588 |
| 0.08 | 0.98196 | 0.98680 | 0.39963 | 0.40445 |
| 0.12 | 1.58572 | 1.59058 | 0.47900 | 0.48382 |
| 0.16 | 2.22574 | 2.23063 | 0.55729 | 0.56211 |
| 0.20 | 2.92767 | 2.93257 | 0.63804 | 0.64286 |
| 0.24 | 3.72458 | 3.72950 | 0.72370 | 0.72852 |
| 0.28 | 4.65765 | 4.66258 | 0.81664 | 0.82146 |
| 0.32 | 5.78589 | 5.79083 | 0.92003 | 0.92485 |
| 0.36 | 7.20334 | 7.20830 | 1.03980 | 1.04462 |
| 0.40 | 9.07731 | 9.08230 | 1.19111 | 1.19593 |

### 3.3 Oscillatory case

Next, we investigate an oscillating $a_t$ case with

$$a_t = A\sin^2\left(\frac{2\pi t}{T}\right) + B, \quad t > 0, \tag{15}$$

where $A > 0$, $T > 0$, and $B > 0$. This kind of $a_t$ is able to model R-tipping in single-species population dynamics in a fluctuating environment (e.g., Abbott et al. [15]). The parameter values are set as $A = 0.8$, $B = 0.01$, and $T = 1$, with which the positivity of $a_t$ is satisfied and it varies in most of the domain $(0, K)$. The other parameter values remain the same as those in the previous computational case.

**Figure 2** shows the computed $X_t$ via the cubature method with $h = 0.0001$ against the initial conditions $X_0 = 0.04i$ ($i = 0,1,2,...,25$). **Figure 2** demonstrates that, unlike the increasing and decreasing $a_t$ cases in **Section 3.2**, there are multiple crossings between the solution and the Allee parameter. The solution trajectories are again separated as those converging to $X_\infty = 1$ and $X_\infty = 0$, and are ordered monotonically, as in (6).

The Euler method with the same computational resolution yields visually comparable solutions as shown in **Figure A1 in Appendix**, and here we compare the Euler and cubature methods in terms of their extinction times $\tau$. **Table 4** compares the values of $\tau$ computed by the two methods for different values of $h$, and **Table 5** summarizes their differences, demonstrating that the Euler method results in earlier extinctions than the cubature method does. Again, the convergence of the computed $\tau$ seems faster in the cubature method, suggesting its greater applicability to the proposed model. The difference between the two methods is more than 20 to 500 times the time increment $h$. This finding combined with the higher accuracy of the cubature method reported in **Appendix** implies that the Euler method may significantly underestimate the extinction time. We consider that this kind of underestimation may also occur in more complex ODEs with non-Lipschitz right-hand sides and caution that the study of extinction time via the Euler method should use a very high computational resolution.

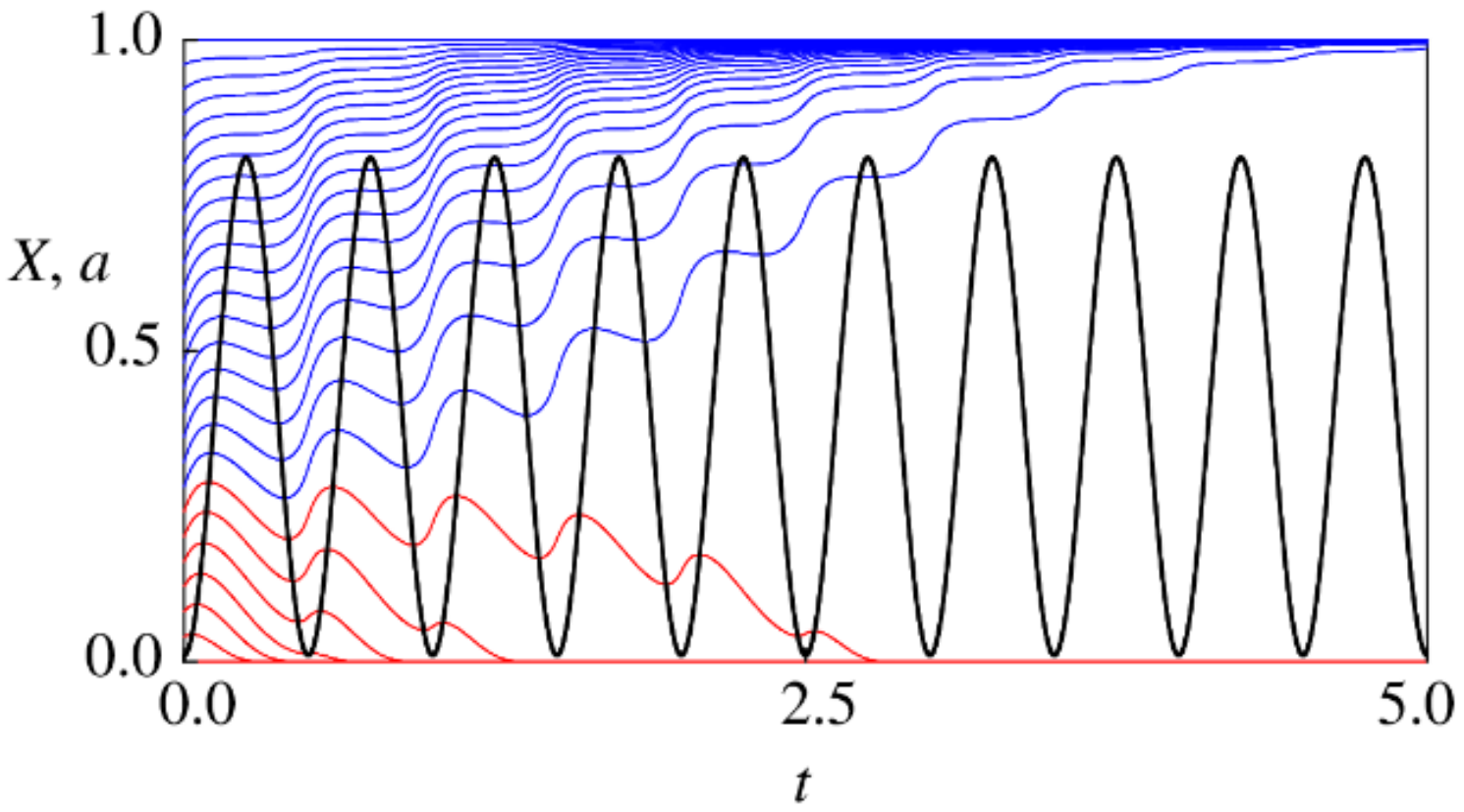

**Figure 2.** The computed $X_t$ by the cubature method. The blue and red curves show the trajectories of $X_t$ with $X_\infty = 1$ (blue) and $X_\infty = 0$ (red). The black curve in each figure panel represents the corresponding trajectory of $a_t$.

**Table 4.** Comparison of the computed $\tau$ with the Euler and cubature methods.

| | $h = 0.001$ | | $h = 0.0001$ | | $h = 0.00001$ | |
|---|---|---|---|---|---|---|
| $X_0$ | Euler | Cubature | Euler | Cubature | Euler | Cubature |
| 0.04 | 0.378 | 0.421 | 0.4058 | 0.4206 | 0.41568 | 0.42050 |
| 0.08 | 0.550 | 0.593 | 0.5771 | 0.5919 | 0.58696 | 0.59177 |
| 0.12 | 0.760 | 0.803 | 0.7867 | 0.8014 | 0.79646 | 0.80127 |
| 0.16 | 1.007 | 1.050 | 1.0336 | 1.0483 | 1.04331 | 1.04811 |
| 0.20 | 1.443 | 1.484 | 1.4663 | 1.4810 | 1.47587 | 1.48067 |
| 0.24 | 2.945 | 2.963 | 2.9260 | 2.9385 | 2.93154 | 2.93612 |

**Table 5.** Difference in the computed $\tau$ between the Euler and cubature methods. The difference here is the computed $\tau$ with the cubature method minus that with the Euler method.

| $X_0$ | $h = 0.001$ | $h = 0.0001$ | $h = 0.00001$ |
|---|---|---|---|
| 0.04 | 0.043 | 0.0148 | 0.00482 |
| 0.08 | 0.043 | 0.0148 | 0.00481 |
| 0.12 | 0.043 | 0.0147 | 0.00481 |
| 0.16 | 0.043 | 0.0147 | 0.00480 |
| 0.20 | 0.041 | 0.0147 | 0.00480 |
| 0.24 | 0.018 | 0.0125 | 0.00458 |

### 3.4 Application to inland fisheries in Japan

We apply the proposed model along with the cubature method to a fisheries-related problem. Inland fisheries are key industries that involve sports and recreational fishing. Habitat degradation of aquatic species is caused by anthropogenetic factors such as water pollution and dam and weir construction and the introduction of nonnative species, and climate change has also been considered the primary factor of the decline in inland fisheries in many countries, such as China [31], Croatia [32], India [33], Southeast Asia [34], and many countries and regions [35]; Japan, which we focus on in this subsection, is not an exception. In Japan, inland fisheries cooperatives serve places to stimulate the economic activities of union members[1]. Inland fisheries cooperatives not only authorize local fishery resources but also manage aquatic environments and play the role of environmental educators, exhibiting multifunctionality toward sustainable coexistence between human and aquatic environments and ecosystems[2]. Union members of each inland fisheries cooperative are allowed to catch specific fish species under certain regulations concerning fishing methods and locations.

Inland fisheries have been degrading in Japan. As reviewed by Nakamura [36], one probable reason for the decline in activity and dissolution of inland fisheries cooperatives in the country is a decrease in membership. We obtained historical data on the total number of union members of fisheries cooperatives in Japan from a public database[3], which can serve as an index to measure the industrial size of inland fisheries in this country. The total number of regular union members was counted every 5 years from 1963 to 2023, as shown in the first and second columns in **Table 6** with a missing value in 2003. The historical data show that the total number of regular union members initially increased since 1960s and then reached a maximum in 1983, after which there was a decreasing trend, which is considered to reflect the decline in activity and dissolution of inland fisheries cooperatives. Importantly, inland fisheries production has decreased since the 1980s in other regions, such as Europe [37], suggesting that the shrinkage of inland fisheries is a global problem. In the case of Japan, the decrease in membership is considered not only due to environmental factors but also to social factors such as the aging of fisheries cooperative staff and union members (e.g., Rahel and Taniguchi [38]; Yoshioka et al. [39]).

We model the total number of regular union members in inland fishery cooperatives in Japan (hereafter called membership for simplicity) by using the ODE (2). Although the membership is an integer and hence discrete variable, it can be approximated as a continuous one if it is significantly larger than 1, which is satisfied in our case, at least during the period when the historical data are available (**Table 6**).

In this subsection, the initial time is set as the beginning of the year 1963, which can be implemented in the proposed model by simply shifting the reference time. The unit of time $t$ is therefore year. We use the following Allee parameter $a_t$, which has preliminarily been found to be suitable for fitting the model against the data:

---

[1] Fisheries Agency in Japan. Sample Articles of Incorporation for Fisheries Cooperatives https://www.jfa.maff.go.jp/j/keiei/gyokyou/mohanteikanrei.html (Last accessed on January 25, 2026. In Japanese)

[2] Fisheries Agency in Japan. Situation surrounding inland fisheries and aquaculture https://www.jfa.maff.go.jp/j/enoki/attach/pdf/naisuimeninfo-45.pdf (Last accessed on January 27, 2026. In Japanese)

[3] E-Stat https://www.e-stat.go.jp/ (Last accessed on January 21, 2026)

$$\ln a_t = \frac{\ln \overline{a} - \ln \underline{a}}{1 + \exp\left(-\frac{t-\theta}{\varepsilon}\right)} + \ln \underline{a}, \quad t \geq 1963, \tag{16}$$

where the meanings of the parameters are the same as those in (14). We use (16) in this subsection to demonstrate that the cubature method works well for sigmoidal functions other than the classical one (14). The Allee parameter $a_t$ quantifies how inland fisheries, such as fishing, recreation, and harvesting are attractive, where a larger value (resp., smaller value) of $a_t$ means that union members become less (resp., more) attracted to these activities.

The parameter values of $X_{1963}$, $K$, $r$, $\ln \overline{a}$, $\ln \underline{a}$, $\varepsilon$, and $\theta$ are estimated by minimizing the following objective function:

$$\frac{1}{M}\sum_{m=1}^{M}\left(X_{\mathrm{Hist},\,t_m} - X_{\mathrm{Theor},\,t_m}\right)^2 \tag{17}$$

by using a common least-squares method, where we impose the following natural constraint to obtain a meaningful model: $\underline{a} \leq \overline{a} \leq K$. Here, the subscripts "Hist" and "Theor" represent historical and theoretical values, respectively, and $M\left(=12\right)$ is the total number of data points $t_m$ ($m = 1,2,3,...,M$) at which the data are available. We compute the theoretical value via the cubature method with $h = 0.001$ (year).

The second through fourth columns of **Table 6** compare historical and theoretical memberships. **Figure 3** shows the historical and theoretical memberships along with the theoretical Allee parameter $a_t$. **Table 7** lists the fitted parameter values. The results obtained show that there is a good agreement between the historical and theoretical results by correctly tracking the rise and fall of the membership, supporting the use of the proposed model in this application problem.

The theoretical Allee parameter $a_t$ implies how the rise and fall of the membership was triggered. The crossing between the membership and the Allee parameter is estimated to be in the late 1983, before and after which there is a rising and falling trend of the historical data, respectively. Hence, the ODE (2) estimates that the maximum in the historical membership time series data is attained near the crossing between the solution and Allee parameter. In the context of our model, this point is understood as the event after which cooperative activities may become less appealing, possibly due to the relative increase of their costs. We find that the time-varying Allee parameter induces a finite-time extinction of the membership in 2051, and that the membership is already close to 0 in the 2040s as shown in **Figure 3**. **Table 7** suggests that the variation in the Allee parameter $a_t$ is extreme in the sense that it nearly varies from 0 to $K$. This finding suggests that the decline in membership was driven by factors that are difficult to avoid. Moreover, the theoretical profile of $a_t$ suggests that the decline of inland fisheries from the standpoint of membership was driven by a gradual process with a transition width of $O\left(10\right)$ years; however, the reason for this gradual variation cannot be narrowed down unless another model(s) for environmental and/or social processes, such as the habitat conditions of aquatic species and nature–human interactions, are explicitly accounted for in the model. Nevertheless, the results obtained would provide some hints for the development of a more sophisticated model. For example, the rise and fall of the annual catch of the fish *Plecoglossus altivelis*[4], a major inland fishery resource in Japan, seem to be correlated with the dynamics of the membership, and can be coupled with the proposed model in future.

Finally, on the basis of a statistical regression analysis, Nakamura [36] estimated that the maximum membership was in 1994 and predicted that extinction occurs in 2035 to 2036; we have not found any other predictions in the literature thus far. The prediction by the proposed model is therefore more conservative, but both the present and previous studies imply that inland fisheries in Japan will be abandoned in the near future with a monotonic shrinkage trend, suggesting that some policy development to promote inland fisheries is necessary so that we can avoid the extinction, e.g., financial incentives for revitalizing activities by inland fisheries cooperatives and their union members to invite young people to participate. The improvement of efficiency to create a comfortable environment for aquatic species including fishery resources, such as feeding grounds and spawning sites, will also be important.

***Remark 2.*** Nakamura [36] also studied the dynamics of membership in inshore fisheries in Japan. We did not analyze the se data because they have been monotonically decreasing since the 1960s, which can be considered a solution always below the Allee parameter in view of the proposed model.

[4] Ministry of Land, Infrastructure, Transport and Tourism, Chubu Regional Development Bureau Water Agency, Chubu Branch https://www.water.go.jp/chubu/nagara/20_followup/pdf/followup_R02_houkokusho.pdf
In Japanese (Last accessed on January 27, 2026)

**Table 6.** Historical and theoretical total number of regular union members. The relative differences are also presented.

| Year | Historical | Theoretical | Relative difference |
|---|---|---|---|
| 1963 | 215,384 | 229,484 | 0.0655 |
| 1968 | 429,258 | 395,669 | 0.0782 |
| 1973 | 487,556 | 503,425 | 0.0325 |
| 1978 | 546,595 | 558,414 | 0.0216 |
| 1983 | 567,113 | 576,248 | 0.0161 |
| 1988 | 563,926 | 568,893 | 0.0088 |
| 1993 | 558,797 | 544,035 | 0.0264 |
| 1998 | 539,089 | 506,384 | 0.0607 |
| 2003 | NA* | 458,802 | NA |
| 2008 | 380,401 | 403,087 | 0.0596 |
| 2013 | 329,239 | 340,561 | 0.0344 |
| 2018 | 271,167 | 272,616 | 0.0053 |
| 2023 | 213,182 | 201,371 | 0.0554 |

* NA: not available.

**Table 7.** Fitted parameter values.

| Parameter | Value |
|---|---|
| $X_{1963}$ (ind) | 229,484 |
| $r$ (1/year) | 0.01691 |
| $K$ (ind) | 1,590,459 |
| $\bar{a}$ (ind) | 1,590,458 |
| $\underline{a}$ (ind) | $3.610 \times 10^{-75}$ * |
| $\varepsilon$ (year) | 10.65 |
| $\theta$ (year) | 1929.17 |

* This small value is considered due to estimating $\ln \underline{a}$ but not $\underline{a}$ itself.

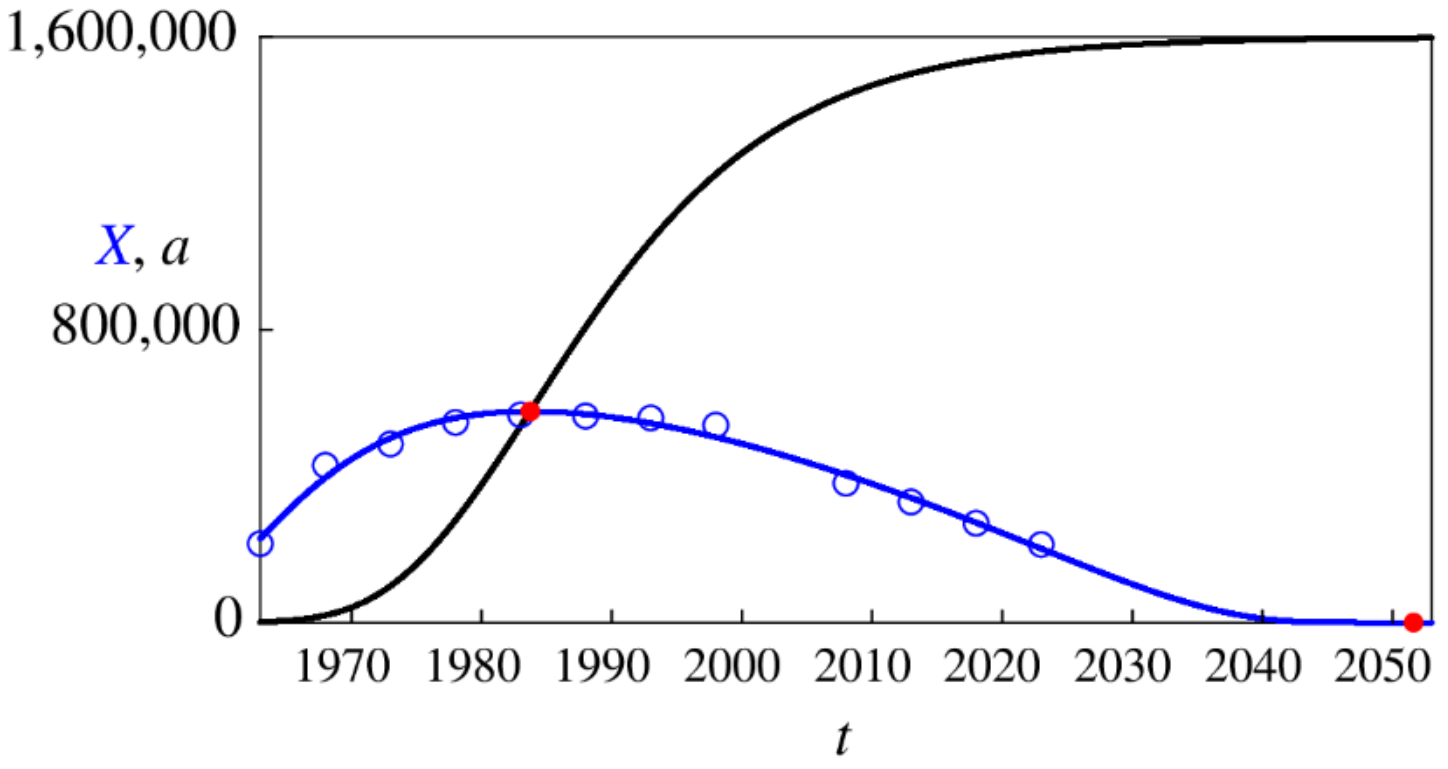

**Figure 3.** Historical (blue blank circles) and theoretical (blue curve) total number of union members $X_t$. The theoretical Allee parameter $a_t$ is also presented in the figure. Their units is "ind". The red filled circles represent the crossing (left) and extinction events (right).

## 4. Conclusion

In this paper, we proposed a simple ODE model that admits an explicit solution. Rete-induced tipping due to time-dependent Allee effects was studied analytically based on the proposed model along with its extinction time. We also discussed a fisheries application of the proposed model. The model was successfully fitted to real data and suggested when and how the rise and fall of inland fisheries in this country occurred. The dynamical system-based analysis based on the proposed model is consistent with the rise and fall of membership in inland fisheries in Japan as a R-tipping phenomenon with a monotone, large-scale, gradual increase in the Allee parameter. This type of analysis connecting R-tipping to industrial problems is still not abundant, and we expect that this study will facilitate such multidisciplinary research.

A limitation of this study is that it only addresses one-dimensional ODE, whereas many multidimensional ODE systems with Allee effects exist. A possible future direction of the proposed model is an extension to a predator–prey model with finite-time extinction. Another direction of the proposed model is to add noise or uncertainty to the right-hand side, which would find relevant applications in biology and ecology where population dynamics in random environments are considered an important research topic. Such a model may be developed through an application of the random variable transformation method [40], with which the tractability of the proposed model can be maintained by choosing a suitable randomization scheme. One may also consider some bounded diffusion processes, such as Jacobi diffusion [41] as a model for the Allee parameter. Piecewise constant noise processes such as the telegraph process can also be considered [42]. In these views, the development of a high-order scheme applicable to the proposed and extended models will be an interesting future topic.

From an application viewpoint, investigations of inland fisheries in countries or regions other than Japan based on the proposed model with suitable extensions would be able to reveal both local and global fisheries trends. Coupling of the proposed model with the resource dynamics mentioned in **Section 3.4** will also be an important research topic, which is currently under investigation by the author.

**• Availability of data and materials**

The data will be made available upon reasonable request to the author.

**• Competing Interests**

The author declares that he has no competing interests.

**• Funding**

This study was supported by the Japan Science and Technology Agency (PRESTO No. JPMJPR24KE) and the Japan Society for the Promotion of Science (KAKENHI No. 25K00240; No. 25K07931).

**• Authors' contributions**

Only the author contributed to this study.

**• Acknowledgements**

None

## Appendix

### A1. Proofs

**Proof of Proposition 1**

First, the unique existence of a local solution that is continuously differentiable near the initial time is guaranteed because the right-hand side of (2) is locally Lipschitz continuous for $X_t > 0$ (Proposition 6.1 in Magnus [43]).

Second, we apply the transformation of variables to the ODE (2). Note that this ODE can be rewritten as the logistic-type ODE

$$\frac{\mathrm{d}(\ln X_t)}{\mathrm{d}t} = r(\ln K - \ln X_t)(\ln X_t - \ln a_t), \quad t > 0. \tag{18}$$

Then, if $X_t \in (0, K)$, then the ODE for the transformed variable $W_t = \dfrac{1}{\ln K - \ln X_t}$ is obtained as

$$\frac{\mathrm{d}W_t}{\mathrm{d}t} = r\left(-1 + \ln\left(\frac{K}{a_t}\right)W_t\right), \quad t > 0 \tag{19}$$

with $W_0 = \dfrac{1}{\ln K - \ln X_0} > 0$. The linear ODE (19) is exactly solved as

$$\begin{aligned} W_t &= W_0 e^{r\int_0^t \ln\left(\frac{K}{a_s}\right)\mathrm{d}s} - r\int_0^t e^{r\int_s^t \ln\left(\frac{K}{a_u}\right)\mathrm{d}u}\mathrm{d}s \\ &= e^{r\int_0^t \ln\left(\frac{K}{a_s}\right)\mathrm{d}s}\left(W_0 - r\int_0^t e^{-r\int_0^s \ln\left(\frac{K}{a_u}\right)\mathrm{d}u}\mathrm{d}s\right), \quad t > 0. \\ &= I_t e^{r\int_0^t \ln\left(\frac{K}{a_s}\right)\mathrm{d}s} \end{aligned} \tag{20}$$

This $W_t$ is positive as long as $I_t > 0$; hence, $X_t > 0$ is guaranteed for $0 \le t < \tau$, where $\tau$ is given by (4) and we can return to $X$ via the inverse transformation $X_t = K\exp(-W_t^{-1})$ and obtain the formula (5) for $0 \le t < \tau$. That for $\tau \ge 0$ follows from

$$\lim_{t \nearrow \tau} X_t = \lim_{t \nearrow \tau} K\exp(-W_t^{-1}) = 0. \tag{21}$$

In addition, $X_t = 0$ and $X_t = K$ are stationary equilibrium of (2).

Finally, the comparison result (6) immediately follows from the representation formula (20) in conjunction with the transformation $X_t = K\exp(-W_t^{-1})$, which completes the proof.

**Proof of Proposition 2**

The proof directly follows due to the boundedness and nonnegativity of the quantity in the exponential function in (11) along with the definition of $k''$ and the formulae (12) and (13).

**A2. Auxiliary computational results**

We present the following auxiliary results:

- ✓ **Table A1** compares the accuracies of the Euler and cubature methods against the extinction time $\tau$ for $a = 0.5$ and $X_0 = 0.32$.
- ✓ **Table A2** compares the accuracies of the two methods against the solution $X$ for $a = 0.5$ and $X_0 = 0.32$.
- ✓ **Table A3** compares the accuracies of the two methods against $\tau$ for $a = 0.5$ and $X_0 = 0.16$.
- ✓ **Table A4** compares the accuracies of the two methods against $X$ for $a = 0.5$ and $X_0 = 0.16$.
- ✓ **Figure A1** shows the $X_t$ computed via the Euler method with the oscillating $a_t$ case for the initial conditions of $X_0 = 0.04i$ ($i = 0,1,2,...,25$).

**Table A1.** Errors of the computed $\tau$ by the Euler and cubature methods for a constant $a\left(=0.5\right)$ case with $X_0 = 0.32$. The convergence rate at the resolution $h$ is estimated as $\log_{10}\left(\frac{\mathrm{Err}_h}{\mathrm{Err}_{h/10}}\right)$, where $\mathrm{Err}_h$ is the absolute error between the exact and computed $\tau$ with the resolution $h$. The exact $\tau$ in this case is 1.35227.

| | $h$ | 0.01 | 0.001 | 0.0001 | 0.00001 |
|---|---|---|---|---|---|
| Error | Euler | 1.023.E-01 | 4.227.E-02 | 1.477.E-02 | 4.804.E-03 |
| | Cubature | 1.773.E-02 | 1.726.E-03 | 1.263.E-04 | 1.630.E-05 |
| Convergence rate | Euler | 3.837.E-01 | 4.566.E-01 | 4.879.E-01 | |
| | Cubature | 1.012.E+00 | 1.136.E+00 | 8.892.E-01 | |

**Table A2.** Average errors of the computed $X$ between the Euler and cubature methods during $0 < t \leq 10$ for a constant $a\left(=0.5\right)$ case with $X_0 = 0.32$. The convergence rate is estimated via the method explained in the caption of **Table A1**.

| | $h$ | 0.01 | 0.001 | 0.0001 | 0.00001 |
|---|---|---|---|---|---|
| Error | Euler | 1.535.E-03 | 1.535.E-04 | 1.534.E-05 | 1.534.E-06 |
| | Cubature | 2.206.E-03 | 2.200.E-04 | 2.199.E-05 | 2.199.E-06 |
| Convergence rate | Euler | 1.000.E+00 | 1.000.E+00 | 1.000.E+00 | |
| | Cubature | 1.001.E+00 | 1.000.E+00 | 1.000.E+00 | |

**Table A3.** The same as **Table A1** but with $X_0 = 0.16$.

| | $h$ | 0.01 | 0.001 | 0.0001 | 0.00001 |
|---|---|---|---|---|---|
| Error | Euler | 1.056.E-01 | 4.256.E-02 | 1.476.E-02 | 4.810.E-03 |
| | Cubature | 4.440.E-03 | 4.404.E-04 | 4.039.E-05 | 1.039.E-05 |
| Convergence rate | Euler | 3.945.E-01 | 4.599.E-01 | 4.870.E-01 | |
| | Cubature | 1.004.E+00 | 1.038.E+00 | 5.897.E-01 | |

**Table A4.** The same as **Table A2** but with $X_0 = 0.16$.

| | $h$ | 0.01 | 0.001 | 0.0001 | 0.00001 |
|---|---|---|---|---|---|
| Error | Euler | 6.297.E-04 | 6.160.E-05 | 6.147.E-06 | 6.146.E-07 |
| | Cubature | 5.442.E-04 | 5.435.E-05 | 5.434.E-06 | 5.434.E-07 |
| Convergence rate | Euler | 1.010.E+00 | 1.001.E+00 | 1.000.E+00 | |
| | Cubature | 1.001.E+00 | 1.000.E+00 | 1.000.E+00 | |

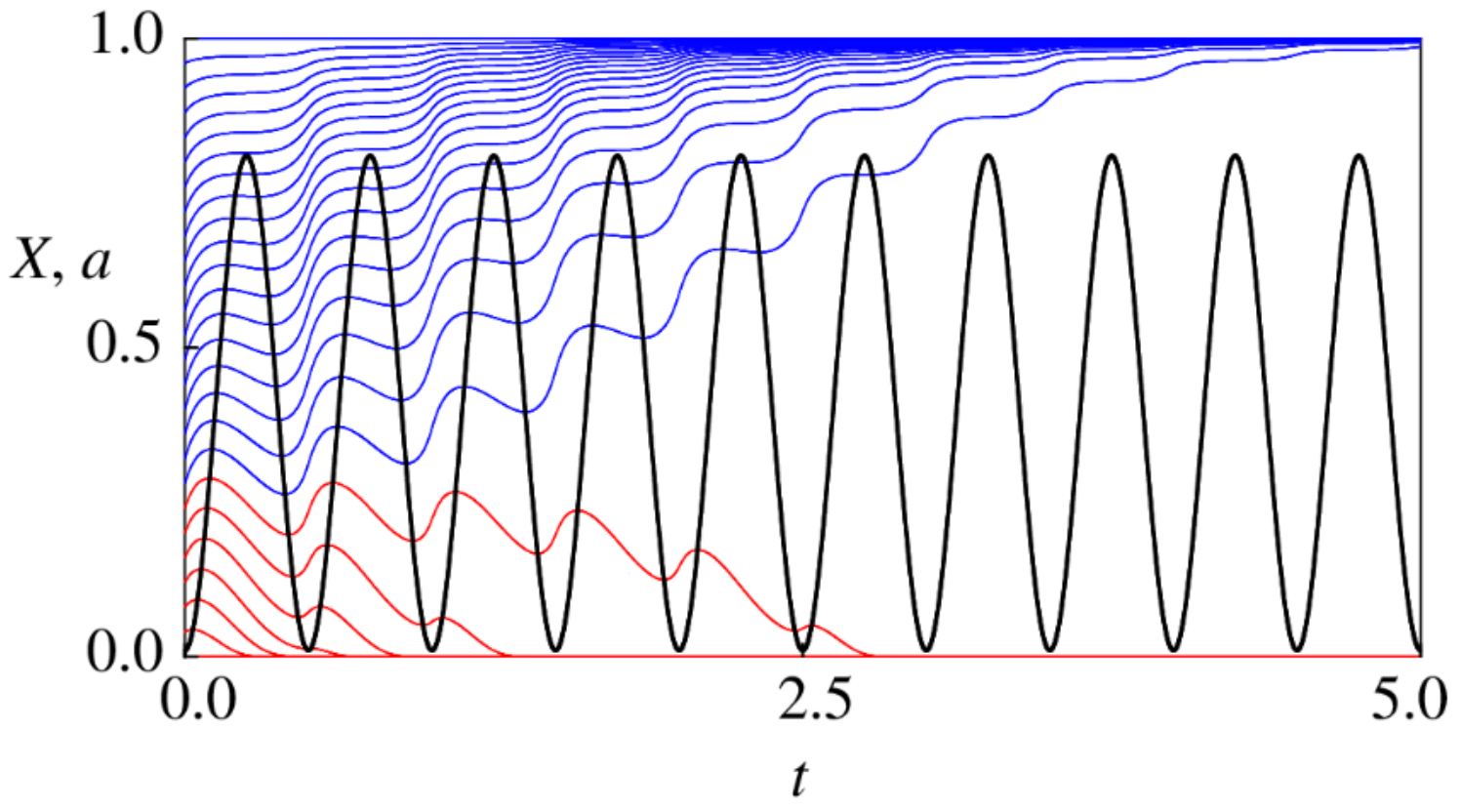

**Figure A1.** The computed $X_t$ by the Euler method for the oscillatory $a_t$ case. The blue and red curves show the trajectories of $X_t$ with $X_\infty = 1$ (blue) and $X_\infty = 0$ (red). The black curve in each figure panel represents the corresponding trajectory of $a_t$.